\documentclass[a4paper,12pt]{amsart}
\usepackage{amssymb}
\usepackage{ifthen}
\usepackage[usenames]{color}
\usepackage{graphicx}
\nonstopmode \numberwithin{equation}{section}
\newtheorem{thm}{Theorem}[section]

\newtheorem{cor}{Corollary}[section]
\newtheorem{lem}{Lemma}[section]
\newtheorem{prop}{Proposition}[section]

\newtheorem{claim}{Claim}[section]
\newtheorem{subclaim}{Subclaim}
\newtheorem{conj}{Conjecture}[section]
\newtheorem{case}{Case}[section]
\newtheorem*{mysolution}{Solution}
\newtheorem{step}{Step}[section]
\theoremstyle{definition}
\newtheorem{defn}{Definition}[section]
\newtheorem{examp}{Example}[section]
\newtheorem{prob}[equation]{Problem}
\newtheorem{ques}{Question}[section]
\newtheorem{rem}{Remark}[section]
\newcounter {own}
\def\theown {\thesection       .\arabic{own}}

\newenvironment{pf}[1][]{%
 \vskip 3mm
 \noindent
 \ifthenelse{\equal{#1}{}}%
  {{\slshape Proof. }}%
  {{\slshape #1.} }%
 }%
{\qed\bigskip}

\newcounter{alphabet}
\newcounter{tmp}
\newenvironment{Thm}[1][]{\refstepcounter{alphabet}%
\bigskip%
\noindent%
{\bf Theorem \Alph{alphabet}}%
\ifthenelse{\equal{#1}{}}{}{ (#1)}%
{\bf .} \itshape}{\vskip 8pt}

\makeatletter

\newcommand{\Ref}[1]{\@ifundefined{r@#1}{}{\setcounter{tmp}{\ref{#1}}\Alph{tmp}}}

\def\be{\begin{equation}}
\def\ee{\end{equation}}

\newcommand{\ben}{\begin{enumerate}}
\newcommand{\een}{\end{enumerate}}

\newcommand{\blem}{\begin{lem}}
\newcommand{\elem}{\end{lem}}
\newcommand{\bthm}{\begin{thm}}
\newcommand{\ethm}{\end{thm}}
\newcommand{\bcor}{\begin{cor}}
\newcommand{\ecor}{\end{cor}}
\newcommand{\beg}{\begin{examp}}
\newcommand{\eeg}{\end{examp}}
\newcommand{\begs}{\begin{examples}}
\newcommand{\eegs}{\end{examples}}
\newcommand{\bdefe}{\begin{defn}}
\newcommand{\edefe}{\end{defn}}
\newcommand{\bprob}{\begin{prob}}
\newcommand{\eprob}{\end{prob}}
\newcommand{\bques}{\begin{ques}}
\newcommand{\eques}{\end{ques}}
\newcommand{\bei}{\begin{itemize}}
\newcommand{\eei}{\end{itemize}}
\newcommand{\bcl}{\begin{claim}}
\newcommand{\ecl}{\end{claim}}
\newcommand{\bscl}{\begin{subclaim}}
\newcommand{\escl}{\end{subclaim}}
\newcommand{\bca}{\begin{case}}
\newcommand{\eca}{\end{case}}
\newcommand{\bstep}{\begin{step}}
\newcommand{\estep}{\end{step}}
\newcommand{\bsol}{\begin{mysolution}}
\newcommand{\esol}{\end{mysolution}}
\newcommand{\bcon}{\begin{conj}}
\newcommand{\econ}{\end{conj}}
\newcommand{\bcons}{\begin{conjs}}
\newcommand{\econs}{\end{conjs}}
\newcommand{\bprop}{\begin{prop}}
\newcommand{\eprop}{\end{prop}}
\newcommand{\br}{\begin{rem}}
\newcommand{\er}{\end{rem}}
\newcommand{\brs}{\begin{rems}}
\newcommand{\ers}{\end{rems}}
\newcommand{\bo}{\begin{obser}}
\newcommand{\eo}{\end{obser}}
\newcommand{\bos}{\begin{obsers}}
\newcommand{\eos}{\end{obsers}}
\newcommand{\bpf}{\begin{pf}}
\newcommand{\epf}{\end{pf}}
\newcommand{\ba}{\begin{array}}
\newcommand{\ea}{\end{array}}
\newcommand{\beq}{\begin{eqnarray}}
\newcommand{\beqq}{\begin{eqnarray*}}
\newcommand{\eeq}{\end{eqnarray}}
\newcommand{\eeqq}{\end{eqnarray*}}

\begin{document}
\bibliographystyle{amsplain}

\title{Dirichlet-type energy of mappings between two concentric annuli}

\author{Jiaolong Chen}
\address{Jiaolong Chen, Key Laboratory of High Performance Computing and Stochastic Information Processing (HPCSIP)
(Ministry of Education of China), School of Mathematics and Statistics, Hunan Normal University, Changsha, Hunan 410081, People's Republic of China}
\email{jiaolongchen@sina.com}

\author{David Kalaj}
\address{David Kalaj, University of Montenegro, Faculty of Natural Sciences and Mathematics, Cetinjski put b.b. 81000 Podgorica, Montenegro}
\email{davidkalaj@gmail.com}


\begin{abstract}
Let $\mathbb{A}$ and $\mathbb{A_{*}}$ be two non-degenerate spherical annuli in $\mathbb{R}^{n}$ equipped with the Euclidean metric
and the weighted metric $|y|^{1-n}$, respectively.
Let $\mathcal{F}(\mathbb{A},\mathbb{A_{*}})$ denote  the class of homeomorphisms in $\mathcal{W}^{1,n-1}(\mathbb{A},\mathbb{A_{*}})$.
For $n=3$, the second author \cite{kalaj2018} proved that the minimizers of the Dirichlet-type energy
$\mathcal{E}[h]=\int_{\mathbb{A}} \frac{\|Dh(x)\|^{n-1}}{|h(x)|^{n-1}}dx$
are certain generalized radial diffeomorphisms,
where $h\in \mathcal{F}(\mathbb{A},\mathbb{A_{*}})$.
For the case $n\geq 4$,
he conjectured that the minimizers are also certain
generalized radial diffeomorphisms between $\mathbb{A}$ and $\mathbb{A_{*}}$.

The main aim of this paper is to consider this conjecture.
First, we investigate the minimality  of the following combined energy integral:
$$
\mathbb{E}[a,b][h]
=\int_{\mathbb{A}}\frac{a^{2}\rho^{n-1}(x)\|DS(x)\|^{n-1}+b^{2}|\nabla \rho(x)|^{n-1}}{|\rho(x)|^{n-1}}dx,
$$
where $h=\rho S\in \mathcal{F}(\mathbb{A},\mathbb{A_{*}})$, $\rho=|h|$ and $a,b>0$.
The obtained result is a generalization of \cite[Theorem 1.1]{kalaj2018}.
As an application,
we show that the above conjecture is almost true for the case $n\geq 4$,
i.e., the minimizer of the energy integral $\mathcal{E}[h]$ does not exist
but there exists
a minimizing sequence which belongs to the generalized radial mappings.
\end{abstract}

\maketitle \pagestyle{myheadings} \markboth{Jiaolong Chen and David Kalaj}{Dirichlet-type energy of mappings between two concentric annuli}

\section{Introduction and statement of the main results}\label{sec-1}
For $n\geq2$, $0< r<R$ and $0<r_{*}<R_{*}$,
let $\mathbb{A}=\{x\in \mathbb{R}^{n}:r<|x|<R\}$ and
$\mathbb{A_{*}}=\{x\in \mathbb{R}^{n}:r_{*}<|x|<R_{*}\}$ be two spherical annuli in $\mathbb{R}^{n}$.
We write $\mathbb{S}^{n-1}(t)=\{x\in \mathbb{R}^{n}:|x|=t\}$ and $\mathbb{S}^{n-1} =\mathbb{S}^{n-1}(1)$.

For natural number $n$, let
$A=(a_{i,j})_{n\times n}\in \mathbb{R}^{n\times  n}$.
We use $A^{T}$ to denote the transpose of $A$.
The {\it Hilbert-Schmit norm}, also called the {\it Frobenius norm},
of $A$ is denoted by $\|A\|$,
where
$$
\|A\|^{2}
=\sum_{1\leq i,j\leq n} \left| a_{i,j} \right|^2
= \text{tr} [A^{T}A].
$$

For $p\geq1$, we say that a mapping $h$ belongs to the class $\mathcal{W}^{1, p}(\mathbb{A},\mathbb{A_{*}})$, if $h$ belongs to the
Sobolev space $\mathcal{W}^{1, p}(\mathbb{A})$ and maps $\mathbb{A}$ onto $\mathbb{A_{*}}$.
Let $h=(h^{1},\ldots,h^{n})$ belong to $\mathcal{W}^{1,p}(\mathbb{A},\mathbb{A_{*}})$.
We denote the {\it Jacobian matrix} of $h$ at the point $x=(x_{1},\ldots,x_{n})$ by $Dh(x)$, where
$Dh(x)=\left(\frac{\partial h^{i}(x)}{\partial x_{j}}\right)_{n\times n}\in \mathbb{R}^{n\times  n}$. Then
$$
\|Dh\|^{2}=\sum_{1\leq i,j\leq n}
\left|  \frac{\partial h^{i}(x)}{\partial x_{j}}\right|^2.
$$
Here $\frac{\partial h^{i}}{\partial x_{j}}$ denotes the weak partial derivatives of $h^{i}$ with respect to $x_{j}$.
If $h$ is continuous and belongs to
$\mathcal{W}^{1, p}(\mathbb{A},\mathbb{A_{*}})$ with $p\geq1$,
then the weak and ordinary partial derivatives coincide a.e. in $\mathbb{A}$ (cf. \cite[Proposition 1.2]{rick}).
Let $h=\rho S$, where $S=\frac{h}{|h|}$ and $\rho=|h|$.
By \cite[Equality (3.2)]{kalaj2018}, we obtain that
$$
Dh(x)=\nabla \rho(x)\otimes S(x) +\rho \cdot DS(x)
$$
and
\begin{equation}\label{eq-1.1}
\|Dh(x)\|^{2}=|\nabla\rho(x)|^{2} +\rho^{2} \|DS(x)\|^{2},
\end{equation}
where $\nabla \rho$ denotes the gradient of $\rho$.

We say that $h:\mathbb{A} \rightarrow\mathbb{A}_{*} $ is a {\it generalized radial mapping} or a {\it quasiradial mapping},
 if there exists a conformal mapping $T$ of $\mathbb{S}^{n-1}$ onto $\mathbb{S}^{n-1}$,
 so that $h(x)=\rho(|x|)T(\frac{x}{|x|})$.
If $T$ is the identity, then we say that $h$ is a {\it radial mapping}.

We use $\mathcal{R}(\mathbb{A},\mathbb{A_{*}})$  to denote the class of radial homeomorphisms in $\mathcal{W}^{1,n-1}(\mathbb{A},\mathbb{A_{*}})$
and use $\mathcal{P}(\mathbb{A},\mathbb{A_{*}})$ to denote the class of generalized radial homeomorphisms in $\mathcal{W}^{1,n-1}(\mathbb{A},\mathbb{A_{*}})$.
We also use $\mathcal{F}(\mathbb{A},\mathbb{A_{*}})$ to
denote the class of homeomorphisms in $\mathcal{W}^{1,n-1}(\mathbb{A},\mathbb{A_{*}})$.

One of interesting and important problems in nonlinear elasticity is whether the radially symmetric
minimizers are indeed global minimizers of the given physically reasonable energy.
For example, Iwaniec and Onninen \cite{iwon} discussed the minimizers of the
following two energy integrals:
$$\mathfrak{E}[h]=\int_{\mathbb{A}} \|Dh(x)\|^{n} dx
\quad\text{and}\quad
 \mathfrak{F}[h]=\int_{\mathbb{A}} \frac{\|Dh(x)\|^{n}}{|h(x)|^{n}} dx
 $$
among all homeomorphisms in $\mathcal{W}^{1,n}(\mathbb{A},\mathbb{A}_{*})$, respectively.
Under some additional conditions,
they proved that the radial minimizers are always the global minimizers.
Later, Koski and Onninen \cite{koski2018} investigated the minimizers of the $p$-harmonic energy
$$\mathfrak{E}_{p}[h]=\int_{\mathbb{X}} \|Dh(x)\|^{p} dx $$
among all homeomorphisms in $\mathcal{W}^{1,p}(\mathbb{X},\mathbb{X}_{*})$,
where $\mathbb{X}$ and $\mathbb{X}_{*}$ are planar annuli and $p\in[1,2)$.
However, if $p=1$ and $\mathbb{X}=\mathbb{X}_{*}=\{x\in \mathbb{R}^2: 0<|x|<1\}$,
Koski and Onninen found that the infimum energy cannot be achieved within the radial mappings.
Recently, the second author \cite{kalaj2018}  studied the Dirichlet-type energy $\mathcal{E}(h)$ among mappings in $\mathcal{F}(\mathbb{A},\mathbb{A_{*}})$,
where
\begin{equation}\label{eq-1.2}
 \mathcal{E}[h]=\int_{\mathbb{A}} \frac{\|Dh(x)\|^{n-1}}{|h(x)|^{n-1}}dx.
\end{equation}
For $n=3$, he proved that the minimizers of $\mathcal{E}(h)$ are certain generalized radial diffeomorphism (cf. \cite[Theorem 1.1]{kalaj2018}).
Motivated by the case $n=3$,  he conjectured
that:

\bcon\label{con-1.1} {\rm (}cf. \cite[Conjecture A.2]{kalaj2018}{\rm )}
For $n\geq 4$, the  energy $\mathcal{E}[h]$ among
mappings in $\mathcal{F}(\mathbb{A},\mathbb{A_{*}})$
achieves its minimum for certain generalized radial diffeomorphisms between $\mathbb{A}$ and $\mathbb{A_{*}}$.
\econ

The main purpose of this paper is to investigate this conjecture.
Recall that if $h(z)=\rho(z)e^{i\Theta(z)}$ is a mapping from planar ring $\mathbb{A}$ onto planar ring $\mathbb{A}_{*}$,
then
$$
\|Dh\|^{2}=|\nabla \rho|^2+\rho^2|\nabla \Theta|^2,
$$
where $\rho=|h|$.
In \cite{ka2020}, the second author studied the minimality  of the combined distortion integral
$$
\mathcal{K}[a,b][h]=\int_{\mathbb{A} }\frac{a^{2} \rho^2(z)|\nabla \Theta(z)|^2+b^{2}|\nabla \rho(z)|^2 }{\det Dh(z)}dz,
$$
where $a$ and $b$ are two positive constants.
Inspired by this paper, first we consider the minimality  of the combined energy integral
\begin{equation}\label{eq-1.3}
\mathbb{E}[a,b][h]
=\int_{\mathbb{A}}\frac{a^{2}\rho^{n-1}(x)\|DS(x)\|^{n-1}+b^{2}|\nabla \rho(x)|^{n-1}}{|\rho(x)|^{n-1}}dx
\end{equation}
among $h=\rho S\in \mathcal{F}(\mathbb{A},\mathbb{A_{*}})$, where  $a,b>0$ and $\rho=|h|$.
Obviously, when $n=3$ and $a=b=1$, it follows from \eqref{eq-1.1} and \eqref{eq-1.2} that $\mathbb{E}[a,b][h]$ coincides with $\mathcal{E}[h]$.

For simplicity,
we define
\begin{equation}\label{eq-1.4}
\;\;
\alpha(t)
=\frac{R^{\frac{1}{n-2} }
 \left(t^{\frac{1}{n-2} }-r^{\frac{1}{n-2} }\right)}
{t^{\frac{1}{n-2} }
\left(R^{\frac{1}{n-2} }-r^{\frac{1}{n-2} }\right)}.
\end{equation}
where $t\in[r,R]$.
We also use
 $h^{\lambda}$ and $h_{i}^{\lambda}$ to denote the mappings in $\mathcal{P}(\mathbb{A},\mathbb{A_{*}})$ with the representations
\begin{equation}\label{eq-1.5}
h^{\lambda}(x)=H(|x|)\Phi^{\lambda}\left(\frac{x}{|x|}\right)
\quad\text{and}\quad
h_{i}^{\lambda}(x)=H_{i}(|x|)\Phi^{\lambda}\left(\frac{x}{|x|}\right),
\end{equation}
respectively,
where $i\in \mathbb{Z}^{+}$, $\lambda>0$ and $\Phi^{\lambda}$ is the mapping from \eqref{eq-2.7} (see below).

The following is our result on $\mathbb{E}[a,b][h]$.

\begin{thm}\label{thm-1.1}
For $n\geq3$ and $h\in\mathcal{F}(\mathbb{A},\mathbb{A_{*}})$,
\begin{equation}\label{eq-1.6}
\begin{split}
\mathbb{E}&[a,b][h]\\
&\geq \omega_{n-1}
\left(a^{2}
(n-1)^{\frac{n-1}{2}}(R-r)+
\frac{b^{2}}{(n-2)^{n-2}}
 \frac{Rr}{ (  R^{\frac{1}{n-2} }-r^{\frac{1}{n-2} })^{n-2}}
 \log^{n-1}\frac{R_{*}}{r_{*}}
\right),
\end{split}
\end{equation}
where $ \omega_{n-1}$ is the $(n-1)$-dimensional
Lebesgue measure of  $\mathbb{S}^{n-1}$.
The equality holds for  the following two generalized radial diffeomorphisms
$$
h_{1}(x)=r_{*}\left( \frac{R_{*}}{r_{*}}\right)^{\alpha(|x|)} T\left(\frac{x}{|x|}\right)
\;
\text{and}
\;\;
h_{2}(x)=R_{*}\left( \frac{r_{*}}{R_{*}}\right)^{\alpha(|x|)}T\left(\frac{x}{|x|}\right),
$$
and the minimizers are unique up to a conformal change $T$
of  $\mathbb{S}^{n-1}$.
\end{thm}

\begin{rem}\label{remark-1.1}
The key tools which we will use in the proof of Theorem \ref{thm-1.1} are the free Lagrangians.
A free Lagrangian is a nonlinear differential $n$-form $G(x,h,Dh)dx$ defined in Sobolev homeomorphisms
$h$ from $\mathbb{A}$ onto $\mathbb{A_{*}} $ whose integral depends only on the homotopy class of $h$ (cf. \cite{iwon}).
\end{rem}

As an application of Theorem \ref{thm-1.1}, we find that Conjecture \ref{con-1.1} is almost true,
i.e., the minimizer of the energy integral $ \mathcal{E}[h]$ does not exist
but there exists
a minimizing sequence which belongs to $\mathcal{P}(\mathbb{A},\mathbb{A}_{*})$.
Our result reads as follows.

\begin{thm}\label{thm-1.2}
For $n \geq4$, we have
\[\begin{split}
\inf_{h\in \mathcal{F}(\mathbb{A},\mathbb{A}_{*})}&
 \mathcal{E}[h]
=\inf_{h\in \mathcal{P}(\mathbb{A},\mathbb{A}_{*})}
 \mathcal{E}[h]
=\lim_{\lambda\rightarrow 0^{+}}
 \mathcal{E}[h_{1}^{\lambda}]
=\lim_{\lambda\rightarrow 0^{+}}
 \mathcal{E}[h_{2}^{\lambda}]\\
&=\omega_{n-1}
\left(
(n-1)^{\frac{n-1}{2}}(R-r)+
\frac{1}{(n-2)^{n-2}}
\frac{Rr}{ (  R^{\frac{1}{n-2} }-r^{\frac{1}{n-2} })^{n-2}}
\log^{n-1}\frac{R_{*}}{r_{*}}
\right),
\end{split}\]
where
\begin{equation}\label{eq-1.7}
h_{1}^{\lambda}(x)
=r_{*}\left( \frac{R_{*}}{r_{*}}\right)^{\alpha(|x|)}\Phi^{\lambda}\left(\frac{x}{|x|}\right)
\;
\text{and}
\;\;
h_{2}^{\lambda}(x)
=R_{*}\left( \frac{r_{*}}{R_{*}}\right)^{\alpha(|x|)}\Phi^{\lambda}\left(\frac{x}{|x|}\right).
\end{equation}
\end{thm}

Further, we establish the following relationships between the Dirichlet-type energy of mappings in $\mathcal{R}(\mathbb{A},\mathbb{A_{*}})$
and  mappings in  $\mathcal{P}(\mathbb{A},\mathbb{A_{*}})$.
Our results are as follows.

\begin{thm}\label{thm-1.3}
For any $n\geq4$ and $h\in  \mathcal{R}(\mathbb{A},\mathbb{A_{*}})$,
there exists a mapping $h^{\lambda}$ such that
 $\mathcal{E}[h^{\lambda}]< \mathcal{E}[h]$ with $|h|=|h^{\lambda}|$ and $\lambda\not=1$.
\end{thm}

\begin{thm}\label{thm-1.4}
For $n\geq4$,
$\inf_{h\in \mathcal{P}(\mathbb{A},\mathbb{A_{*}})}\mathcal{E}[h]
<\min_{h\in \mathcal{R}(\mathbb{A},\mathbb{A_{*}})}\mathcal{E}[h].$
\end{thm}

The rest of this paper is organized as follows.
In Section \ref{sec-2}, some necessary terminology and preliminary results are given.
In Section \ref{sec-3}, Theorem \ref{thm-1.1} is proved, and Theorem \ref{thm-1.2} is shown in Section \ref{sec-4}.
Section \ref{sec-5} is devoted to the proofs of
Theorems \ref{thm-1.3} and \ref{thm-1.4}.
\section{Preliminary results}\label{sec-2}

For any fixed $x\in \mathbb{A}$, we assume that $\{U_{1},\ldots,U_{n-1},N\}$
is a system of mutually orthogonal vectors of the unit norm, where $N=\frac{x}{|x|}$ and
the vectors $\{U_{1},\ldots,U_{n-1}\}$ are arbitrarily chosen.
For $h\in\mathcal{W}^{1, p}(\mathbb{A},\mathbb{A_{*}})$ with $p\geq1$,
define
 $h_{N}(x)=Dh(x)N$ and set $h_{U_{i}}(x)=Dh(x)U_{i}$ for $i\in\{1,\ldots,n-1\}$.
Since the Hilbert-Schimdt norm of $Dh$ is independent of basis
(cf. \cite[Page 8]{cuneo}),
we have
\begin{equation}\label{eq-2.1}
\|Dh\|^{2}
=|h_{N}|^{2}+|h_{U_{1}}|^{2}+|h_{U_{2}}|^{2}+\cdots+|h_{U_{n-1}}|^{2}.
\end{equation}
Let $S=\frac{h}{|h|}$ and define the {\it Gram determinant} of $S$ at $x\in \mathbb{A}$ by
\begin{equation}\label{eq-2.2}
 D_{S}(x)=\left|S_{U_{1}}(x)\times S_{U_{2}}(x)\times \cdots \times S_{U_{n-1}}(x) \right|
\end{equation}
(cf. \cite[Section 3]{kalaj2018}).
Then
\begin{equation}\label{eq-2.3}
\begin{split}\|DS(x)\|^{n-1}
 =&\big(|S_{U_{1}}(x)|^{2}+ |S_{U_{2}}(x)|^{2}+\cdots+|S_{U_{n-1}}(x)|^{2}+|S_{N}(x)|^{2}\big)^{\frac{n-1}{2}}\\ \geq&\big(|S_{U_{1}}(x)|^{2}+ |S_{U_{2}}(x)|^{2}+\cdots+|S_{U_{n-1}}(x)|^{2}\big)^{\frac{n-1}{2}}\\
 \geq& (n-1)^{\frac{n-1}{2}}|S_{U_{1}}(x)|\cdot|S_{U_{2}}(x)|\cdots|S_{U_{n-1}}(x)|\\ \geq& (n-1)^{\frac{n-1}{2}}D_{S}(x).
\end{split}
\end{equation}
Further, we have the following result.

 \begin{Thm}\label{Thm-A} \rm{(}\cite[Lemma 3.1]{kalaj2018})
Let $h\in \mathcal{F}(\mathbb{A},\mathbb{A}_{*})$  and $S=\frac{h}{|h|}$.
Then
$$
\int_{\mathbb{A} } D_{S}(x) dx
\geq  (R-r)\omega_{n-1}.
$$
\end{Thm}

Let $h(x)=\rho(t)S(\xi)$ be a generalized radial mapping in $\mathbb{A}$, where $x=r\xi$, $t\in(r,R)$ and $\xi\in \mathbb{S}^{n-1}$.
It follows from \cite[Page 10]{cuneo} that
$$
h_{N}(x)=\dot{\rho}(r)S(\xi)
\quad\text{and}\quad
h_{U_{i}}(x)= \frac{\rho(r)}{r}S_{U_{i}}(\xi),
$$
where $\dot{\rho} (r)=\frac{d \rho(r)}{d r}$ and
$
 S_{U_{i}}(\xi)=DS(\xi)U_{i}.
$
Hence, we obtain from \eqref{eq-2.1} that
\begin{equation}\label{eq-2.4}
\|Dh(x)\|^{2}
=\dot{\rho}^{2}(r)+\frac{\rho^{2}(r)}{r^{2}}\big(\|DS(\xi)\|^{2}-|S_{N}(\xi)|^{2}\big).
\end{equation}
Note that $S$ is a conformal mapping from $\mathbb{S}^{n-1}$ onto $\mathbb{S}^{n-1}$,
which means that
\begin{equation}\label{eq-2.5}
|S_{N}(\xi)|=0,
\;\;
|S_{U_{1}}(\xi)|=|S_{U_{2}}(\xi)|=\cdots =|S_{U_{n-1}}(\xi)|,
\end{equation}
and
$S_{U_{1}}(\xi),S_{U_{2}}(\xi),\ldots,S_{U_{n-1}}(\xi)$ are mutually orthogonal vectors.
By \eqref{eq-2.2} and \eqref{eq-2.3}, we see that
$$
 D_{ S} (\xi)
=\left|S_{U_{1}}(\xi)\times S_{U_{2}}(\xi)\times \cdots \times S_{U_{n-1}}(\xi) \right|
=|S_{U_{1}}(\xi)|\cdot|S_{U_{2}}(\xi)|\cdots|S_{U_{n-1}}(\xi)|
$$
and
\begin{equation}\label{eq-2.6}
 \|DS(\xi)\|^{n-1}
=\big(|S_{U_{1}}(\xi)|^{2}+ |S_{U_{2}}(\xi)|^{2}+\cdots+|S_{U_{n-1}}(\xi)|^{2} \big)^{\frac{n-1}{2}}
=  (n-1)^{\frac{n-1}{2}} D_{ S}(\xi).
\end{equation}

Let $\Pi:\mathbb{S}^{n-1}\rightarrow\widehat{\mathbb{R}}^{n-1}$ denote the stereographic projection of $\mathbb{S}^{n-1} $ through the south pole onto $\widehat{\mathbb{R}}^{n-1}$,
where $\widehat{\mathbb{R}}^{n-1}= \mathbb{R}^{n-1} \cup \{\infty\}$.
For any $\lambda>0$ and $x\in \widehat{\mathbb{R}}^{n-1}$, let $g_{\lambda}(x)=\lambda x$.
Then it follows from \cite[Proof of Theorem 3.1]{kalaj2019} that
\begin{equation}\label{eq-2.7}
\Phi^{\lambda}=\Pi^{-1}\circ g_{\lambda} \circ\Pi
\end{equation}
 is a conformal mapping  from $\mathbb{S}^{n-1}$ onto $\mathbb{S}^{n-1}$.
 Obviously, if $\lambda=1$, then $\Phi^{\lambda}$ coincides with the identity mapping.
Let $\xi=(\cos \theta,\mathfrak{s} \sin\theta)\in \mathbb{S}^{n-1}$ be a point of longitude $\mathfrak{s} \in \mathbb{S}^{n-2} $
and meridian $\theta\in[0,\pi]$.
The south pole corresponds to $\theta=\pi$.
Then
$$
\Phi^{\lambda}(\xi)=\big(\cos \varphi(\theta),\mathfrak{s} \sin\varphi(\theta)\big)
$$
(cf. \cite[Section 14.3]{iwon}), where
$$
\varphi(\theta)=2\arctan\left(\lambda\tan\frac{\theta}{2}\right).
$$

Furthermore, by \eqref{eq-2.6} and \cite[Equalities (14.38)$\sim$(14.40)]{iwon},
we get that
\begin{equation}\label{eq-2.8}
\|D\Phi^{\lambda}(\xi)\|^{2}
=(n-1)\big( D_{ \Phi^{\lambda}}(\xi)\big)^{\frac{2}{n-1}}
=(n-1)\frac{\sin^{2}\varphi}{\sin^{2}\theta}
\end{equation}
and
\begin{equation}\label{eq-2.9}
\int_{\mathbb{S}^{n-1}} \|D\Phi^{\lambda}(\xi)\|^{n-1}d\sigma(\xi)=(n-1)^{\frac{n-1}{2}} \omega_{n-1},
\end{equation}
where
 $\sigma$ denotes the $(n-1)$-dimensional Lebesgue measure so that $\sigma(\mathbb{S}^{n-1})=\omega_{n-1}$.

\medskip

\section{Proof of Theorem \ref{thm-1.1}}\label{sec-3}
The aim of this section is to prove Theorem \ref{thm-1.1}.
For  convenience, in the rest of this paper, we set $E=(r,R)$ and $F=(r_{*},R_{*})$.

\subsection{Proof of Theorem \ref{thm-1.1} }
Let $h=\rho S\in \mathcal{F}(\mathbb{A},\mathbb{A_{*}})$, where $\rho=|h|$ and $S=\frac{h}{|h|}$.
Before we go to the detailed proof, let us make one shortcut. For every constant $a,b,c>0$,
we have
\begin{equation}\label{eq-3.1}
\mathbb{E}[a,b]\left[\frac{ch}{|h|^2}\right]
=\mathbb{E}[a,b][h]
=\int_{\mathbb{A}} \left(a^{2}
\|D S(x)\|^{n-1}+b^{2}\frac{|\nabla\rho(x)|^{n-1}}{\rho^{n-1}(x)} \right)dx.
\end{equation}

In order to prove this statement, we let $f=\frac{ch}{|h|^2}$.
Then
$$
|f(x)|=\frac{c}{\rho(x)}
\quad\text{and}\quad
\frac{f(x)}{|f(x)|}
= \frac{h(x)}{|h(x)|}=S(x).
$$
By calculations, we get that
$\nabla |f(x)|=-c\frac{\nabla \rho(x)}{\rho^{2}(x) }$,
and so,
$$
\frac{\big|\nabla |f(x)|\big| }{ |f(x)| }
=\frac{ |\nabla \rho (x)  | }{\rho (x)}.
$$
Then the above equalities and \eqref{eq-1.3} imply that
\eqref{eq-3.1} holds true.
Thus, in the following, we can assume that
\begin{equation}\label{eq-3.2}
\lim_{|x|\rightarrow r}
|h(x)|=r_{*}\;\;\;\text{and}\;\;\;\lim_{|x|\rightarrow R}|h(x)|=R_{*}.
\end{equation}

In order to find the minimizer of $\mathbb{E}[h]$, first, we estimate the integral $\int_{\mathbb{A}}\|DS(x)\|^{n-1}dx$.
It follows from \eqref{eq-2.3} and Theorem \Ref{Thm-A}  that
\begin{equation}\label{eq-3.3}
\int_{\mathbb{A}}  \|D S(x)\|^{n-1}dx
\geq (n-1)^{\frac{n-1}{2}}\int_{\mathbb{A}}D_{S}(x)dx
\geq(R-r)(n-1)^{\frac{n-1}{2}} \omega_{n-1}.
\end{equation}

Second, we estimate the integral $\int_{\mathbb{A}}\frac{|\nabla\rho(x)|^{n-1}}{\rho^{n-1}(x)}dx$.
Let $x=t\xi$, where $t=|x|$ and $\xi\in \mathbb{S}^{n-1}$.
Then
\begin{equation}\label{eq-3.4}
 |\rho_{N}(x) |=
\big|\langle \nabla\rho(x),\xi\rangle\big|
\leq  |\nabla\rho(x) |,
\end{equation}
where $\rho_{N}$ denotes the differentiation of $\rho$ in the direction $N$
and
$$
\rho_{N}(x)=\langle \nabla\rho(x),N\rangle
=\frac{\partial\rho(x)}{\partial t}.
$$
By H\"{o}lder's inequality and \eqref{eq-3.4}, we see that
\begin{equation}\label{eq-3.5}
\int_{\mathbb{A}}\frac{|\nabla\rho(x)|^{n-1}}{\rho^{n-1}(x)}dx
 \geq \left(\int_{\mathbb{A}} \frac{ |\rho_{N}(x) |}{\rho(x)|x|^{n-1}} dx\right)^{n-1}
\left(\int_{\mathbb{A}} \frac{1}{|x|^{\frac{(n-1)^2}{n-2}} } dx\right)^{2-n}.
\end{equation}
The equality holds if and only if
$$
\frac{ |\nabla\rho(x)|}{\rho(x) }
=\frac{ |\rho_{N}(x) |}{\rho(x) }
=C_{1,1} |x|^{-\frac{n-1}{n-2}}
$$
a.e. in $\mathbb{A}$ for some constant $C_{1,1}>0$
(cf. \cite[Page 6]{kuang}).
 Since $h =\rho S$
is a homeomorphism in $\mathcal{W}^{1, n-1}(\mathbb{A},\mathbb{A}_{*})$,
we obtain that (cf. \cite[Equality (2.6)]{ka2020})
\begin{equation}\label{eq-3.6}
\begin{split}
\int_{\mathbb{A}} \frac{ |\rho_{N}(x) |}{\rho(x)|x|^{n-1}} dx
 \geq&\int_{\mathbb{A}} \frac{  \rho_{N}(x) }{\rho(x)|x|^{n-1}} dx
= \int_{\mathbb{S}^{n-1}}\int_{r}^{R} \frac{  \frac{\partial}{\partial t} \rho(t\xi) }{\rho(t\xi) } dtd\sigma(\xi)\\
 =&\omega_{n-1}\int_{r_{*}}^{R_{*}}  \frac{  d\rho }{\rho} =\omega_{n-1}\log \frac{R_*}{r_*}.
\end{split}
\end{equation}
Combing \eqref{eq-3.1}, \eqref{eq-3.3}, \eqref{eq-3.5} and \eqref{eq-3.6}, we see that \eqref{eq-1.6} holds true.

To prove the equality statement, assume that the equalities are attained in all inequalities.
If the equality is attained in \eqref{eq-3.4}, then
$$
|\rho_{N}(x)|
=\left|\left\langle \nabla\rho(x), \frac{x}{|x|}\right\rangle\right|
= |\nabla\rho(x)|,
$$
which implies that the directional derivative of $ \rho_{U_{i}} $ are all vanished since
$$
| \rho_{U_{i}}(x)|=|\langle \nabla\rho(x),U_{i} \rangle|=0
$$
for $i=1,\ldots,n-1$.
Hence, $\rho(x)$ depends only on the radial part of $x=t\xi$, i.e., $\rho$ can be expressed as $\rho(x)=H(t)$,
where $t=|x|$ (cf. \cite[Page 976]{iwon2009}).
If the equalities are attained in \eqref{eq-3.3}, then
we see from \eqref{eq-2.3} that
\begin{center}
$|D_{N}S(x)|=0$, $|D_{U_{1}}S(x)|=|D_{U_{2}}S(x)|=\cdots =|D_{U_{n-1}}S(x)|$ \end{center}
and
$D_{U_{1}}S(x), D_{U_{2}}S(x),\cdots D_{U_{n-1}}S(x)$ are mutually orthogonal vectors.
Thus there exists a conformal mapping $T$ from $\mathbb{S}^{n-1}$ onto $\mathbb{S}^{n-1}$ such that $S(t\xi)=T(\xi)$, where $\xi\in \mathbb{S}^{n-1}$.
Therefore, when the equalities in \eqref{eq-3.3} and \eqref{eq-3.4} hold,
$h$ must be of the form
$$
h(t\xi)=H(t)T(\xi),
$$
where  $t\in E$ and $\xi\in \mathbb{S}^{n-1}$.
This means that $H$ is a homeomorphism  from $E$ onto $F$ in $\mathcal{W}^{1, n-1}(E)$.
Recall that \eqref{eq-3.5} and \eqref{eq-3.6} are equalities if and only if
$$
\frac{\dot{H}(t) }{H(t)}= C_{1,1}  t^{-\frac{n-1}{n-2}}
 $$
a.e. in $E$ for some $C_{1,1}>0$.
Let
 $$
 G(t)=\log H(t)+(n-2)C_{1,1}  t^{-\frac{1}{n-2}}
 $$
 be a continuous mapping in $E$.
Then $G\in \mathcal{W}^{1, n-1}(E)$ and $\dot{G}=0$ a.e. in $E$.
Thus there exists a constant $C_{1,2}$ such that $G=C_{1,2}$ a.e. in $E$.
Since $G$ is continuous in $E$, we obtain that $G\equiv C_{1,2}$ in $E$, i.e.,
$$
H(t)=C_{1,3}\exp\left\{C_{1,4}t^{\frac{1}{2-n}}\right\}
$$
for some constants $C_{1,3}>0$ and $C_{1,4}\in \mathbb{R}$.
This, together with \eqref{eq-3.2}, implies that
$$
C_{1,3}
=r_{*}\left( \frac{R_{*}}{r_{*}}\right)^{\frac{  R^{\frac{1}{n-2} }  }
{ R^{\frac{1}{n-2} }-r^{\frac{1}{n-2} } }}
\quad\text{and}\quad
C_{1,4}
= \frac{(R r)^{\frac{1}{n-2}} }{  R^{\frac{1}{n-2}}- r^{\frac{1}{n-2}} }
\log\frac{r_{*}}{R_{*}}.
$$
Thus the proof of the theorem is complete.
\qed

\section{Proof of Theorem \ref{thm-1.2} }\label{sec-4}
We shall prove Theorem  \ref{thm-1.2} in this section.
The proof will be based on three lemmas.
The first lemma reads as follows.

\begin{lem}\label{lem-4.1}
For any $a\geq b\geq0$ and $s\geq1$, we have
$$
a^{s}-b^{s}\leq s(a-b)(a^{s-1}+b^{s-1}).
$$
\end{lem}
\bpf
If $b=0$, obviously, we have $a^{s}-b^{s}\leq s(a-b)(a^{s-1}+b^{s-1})$.
In the following, we assume that $b>0$.
Note that
$$
s(a-b)(a^{s-1}+b^{s-1}) =s(a^{s}+ab^{s-1}-ba^{s-1}-b^{s}).
$$
Hence, it suffices to prove
\begin{equation}\label{eq-4.1}
\left(\frac{a}{b}\right)^{s}-1
\leq s\left(
\left(\frac{a}{b}\right)^{s}+ \frac{a}{b}-\left(\frac{a}{b}\right)^{s-1}-1
\right).
\end{equation}

Let
$$
g(x)= s\left(x^{s}+ x-x^{s-1}-1\right)-x^{s}+1,
$$
where $x\geq1$.
Obviously,
$$g'(x)= s(s-1)x^{s-1}-s(s-1)x^{s-2}+s$$
and
$$g''(x)= s(s-1)x^{s-3}\big((s-1)x-(s-2)\big).$$
Since $g''(x)\geq0$ in $[1,+\infty)$ and $g'(1)=s>0$,
we see that $g$ is increasing in $[1,+\infty)$.
It follows from $g(1)=0$ that $g(x)\geq0$ in $[1,+\infty)$.
Therefore, \eqref{eq-4.1} holds true,
and the proof of the lemma is complete.
\epf

Next, we establish a general integral representation formula for $\mathcal{E}[h^{\lambda}]$, where $h^{\lambda}$ is the mapping from \eqref{eq-1.5}.

\begin{lem}\label{lem-4.2}
For $n \geq4$ and $\lambda>0$,
\begin{equation}\label{eq-4.2}
\begin{split}
&\mathcal{E}[h^{\lambda}]\\
 &=2^{n-1}\omega_{n-2}\int_{r}^{R} \int_{0}^{+\infty}\left( \frac{t^{2}\dot{H}^{2}(t)}{H^{2}(t)}+
 (n-1)\frac{\lambda^{2} (1+y^{2})^{2} }{(1+\lambda^{2} y^{2})^{2}} \right)^{\frac{n-1}{2}} \frac{y^{n-2}}{(1+y^{2})^{n-1}}dy dt.
\end{split}
\end{equation}
In particular,
\begin{equation}\label{eq-4.3}
 \mathcal{E}[h^{1}]
=\omega_{n-1}\int_{r}^{R}\left(n-1+\frac{t^{2}\dot{H}^{2}(t)}{H^{2}(t)}\right)^{\frac{n-1}{2}} dt .
\end{equation}
\end{lem}
\bpf
For $\lambda>0$,
by \eqref{eq-1.2}, \eqref{eq-1.5}, \eqref{eq-2.4} and \eqref{eq-2.5},
we obtain
\begin{equation}\label{eq-4.4}
\mathcal{E}[h^{\lambda}]
=\int_{r}^{R} \int_{\mathbb{S}^{n-1}}\left( \frac{t^{2}\dot{H}^{2}(t)}{H^{2}(t)}+
\|D\Phi^{\lambda}(\xi)\|^{2} \right)^{\frac{n-1}{2}}d\sigma(\xi)dt.
\end{equation}
By \eqref{eq-2.8} and \eqref{eq-4.4}, we get
$$
\mathcal{E}[h^{\lambda}]
=\omega_{n-2}\int_{r}^{R} \int_{0}^{\pi}\left( \frac{t^{2}\dot{H}^{2}(t)}{H^{2}(t)}+
(n-1)\frac{\sin^{2}\big(2\arctan(\lambda\tan\frac{\theta}{2})\big)}{\sin^{2}\theta} \right)^{\frac{n-1}{2}} \sin^{n-2}\theta \;d\theta dt.
$$
Let $y=\tan\frac{\theta}{2}$ and $\vartheta=\arctan(\lambda y)$.
Then
$\sin\theta=\frac{2y}{1+y^{2}}$,
$d\theta=\frac{2dy}{1+y^{2}}$
and
$$
\sin\left(2\arctan\left(\lambda\tan\frac{\theta}{2}\right)\right)
=\sin \big(2\arctan(\lambda y)\big)
=\frac{2\lambda y}{1+\lambda^{2} y^{2}}.
$$
Combining the above equalities, we see that \eqref{eq-4.2} holds true.
Further,
\begin{equation}\label{eq-4.5}
\begin{split}
\int_{0}^{+\infty} &\frac{y^{n-2}}{(1+y^{2})^{n-1}}dy
=\frac{1}{2}\int_{0}^{+\infty} \frac{y^{\frac{n-3}{2}}}{(1+y)^{n-1}}dy
=\frac{1}{2}\int_{1}^{+\infty}  \frac{(s-1)^{\frac{n-3}{2}}}{s^{n-1}} ds\\
 =&\frac{1}{2}\int_{0}^{1} s^{\frac{n-3}{2}} (1-s)^{\frac{n-3}{2}}ds
=\frac{\Gamma^{2}(\frac{n-1}{2})}{2\Gamma(n-1)}
=\frac{ \sqrt{\pi}\Gamma(\frac{n-1}{2})}{2^{n-1}\Gamma(\frac{n}{2})}
\end{split}
\end{equation}
and
\begin{equation}\label{eq-4.6}
\frac{\omega_{n-2}}{\omega_{n-1}}
=\frac{1}{ \int_{0}^{\pi} \sin^{n-2}\theta d\theta }
=\frac{\Gamma(\frac{n}{2})}{\sqrt{\pi}\Gamma(\frac{n-1}{2})}.
\end{equation}
Then \eqref{eq-4.2}, \eqref{eq-4.5} and \eqref{eq-4.6} yield that \eqref{eq-4.3} holds true.
The proof of the lemma is complete.
\epf

Based on Lemmas \ref{lem-4.1} and \ref{lem-4.2}, we have the following result.

\begin{lem}\label{lem-4.3}
For $n \geq4$, we have
 $$
 \inf_{\lambda\in(0,+\infty)}
 \mathcal{E}[h^{\lambda}]
 = \lim_{\lambda\rightarrow 0^{+}}
 \mathcal{E}[h^{\lambda}]
  <\mathcal{E}[h^{1}]
 $$
and
\[\begin{split}
\inf_{\lambda\in(0,+\infty)} \mathcal{E}[h^{\lambda}]
 \geq&\omega_{n-1}
\left((n-1)^{\frac{n-1}{2}} +\frac{\log^{n-1}\frac{R_{*}}{r_{*}}}{(n-2)^{n-2}}
\frac{Rr}{ (  R^{\frac{1}{n-2} }
-r^{\frac{1}{n-2} })^{n-2}}\right)\\
 =& \lim_{\lambda\rightarrow 0^{+}}
 \mathcal{E}[h_{1}^{\lambda}]
=\lim_{\lambda\rightarrow 0^{+}}
 \mathcal{E}[h_{2}^{\lambda}],
\end{split}\]
where $h_{1}^{\lambda}$
and $h_{2}^{\lambda}$ are the mappings from
 \eqref{eq-1.7}.
\end{lem}
\bpf
In order to prove the lemma, first, we show the following claim.
\begin{claim}\label{claim-4.1}
For $n \geq4$,
\begin{eqnarray*}
 \inf_{\lambda\in(0,+\infty)}
 \mathcal{E}[h^{\lambda}]
  = \lim_{\lambda\rightarrow 0^{+}}
 \mathcal{E}[h^{\lambda}]
=
\omega_{n-1} \int_{r}^{R}
 \left( (n-1)^{\frac{n-1}{2}} +\frac{t^{n-1}\dot{H}^{n-1}(t)}{H^{n-1}(t)}\right)  dt
< \mathcal{E}[h^{1}].
\end{eqnarray*}
\end{claim}
For any $\lambda>0$, by calculations, we get
\begin{eqnarray*}
\quad\int_{r}^{R}\int_{0}^{+\infty} \frac{\lambda^{n-1} (1+y^{2})^{n-1} }
{(1+\lambda^{2} y^{2})^{n-1}}  \frac{y^{n-2}}{(1+y^{2})^{n-1}}dydt
=\int_{r}^{R}\int_{0}^{+\infty}  \frac{s^{n-2}}{(1+s^{2})^{n-1}}dsdt.
\end{eqnarray*}
This, together with
\eqref{eq-1.5}, \eqref{eq-4.2}, \eqref{eq-4.5} and \eqref{eq-4.6}, implies that
\begin{equation}\label{eq-4.7}
\begin{split}
\mathcal{E} [h^{\lambda}]
\geq&2^{n-1}\omega_{n-2}\int_{r}^{R}\int_{0}^{+\infty}
\left( \frac{t^{n-1}\dot{H}^{n-1}(t)}{H^{n-1}(t)}+ (n-1)^{\frac{n-1}{2}}
\frac{\lambda^{n-1} (1+y^{2})^{n-1} }{(1+\lambda^{2} y^{2})^{n-1}} \right)\\ &\times \frac{y^{n-2}}{(1+y^{2})^{n-1}}dydt\\
=& \omega_{n-1} \int_{r}^{R} \left( (n-1)^{\frac{n-1}{2}} + \frac{t^{n-1}\dot{H}^{n-1}(t)}{H^{n-1}(t)}\right)  dt.
\end{split}
\end{equation}

In the following, we prove that
\begin{equation}\label{eq-4.8}
\lim_{\lambda\rightarrow 0^{+}}
 \mathcal{E}[h^{\lambda}]
=\omega_{n-1} \int_{r}^{R} \left((n-1)^{\frac{n-1}{2}}+ \frac{t^{n-1}\dot{H}^{n-1}(t)}{H^{n-1}(t)}\right)  dt.
\end{equation}
By \eqref{eq-4.2} and \eqref{eq-4.7}, we see that \eqref{eq-4.8} is equivalent to
\begin{equation}\label{eq-4.9}
\begin{split}
\lim_{\lambda\rightarrow 0^{+}}& \int_{r}^{R} \int_{0}^{+\infty}
\bigg[\left( \frac{t^{2}\dot{H}^{2}(t)}{H^{2}(t)}+
 (n-1)\frac{\lambda^{2} (1+y^{2})^{2} }{(1+\lambda^{2} y^{2})^{2}} \right)^{\frac{n-1}{2}} -\frac{t^{n-1}\dot{H}^{n-1}(t)}{H^{n-1}(t)} \\
 & -  (n-1)^{\frac{n-1}{2}}\frac{\lambda^{n-1} (1+y^{2})^{n-1} }{(1+\lambda^{2} y^{2})^{n-1}}
\bigg] \frac{y^{n-2}}{(1+y^{2})^{n-1}}dy dt=0.
\end{split}
\end{equation}

For fixed $a\geq 0$,
obviously,  the mapping
$$x\mapsto (a^{2}+x^{2})^{\frac{n-1}{2}}-a^{n-1}-x^{n-1}$$
 is increasing  in $[0,+\infty)$
and
$$
\frac{\lambda^{2}  }{(1+\lambda^{2} y^{2})^{2}} \leq \frac{1}{4y^{2}}
$$
for any $y>0$ and $\lambda\geq0$.
 Then for any $t\in E$, $y>0$ and $\lambda\geq0$,
 it follows from Lemma \ref{lem-4.1} that
\[\begin{split}
\bigg( &\frac{t^{2}\dot{H}^{2}(t)}{H^{2}(t)}+
 (n-1)\frac{\lambda^{2} (1+y^{2})^{2} }{(1+\lambda^{2} y^{2})^{2}} \bigg)^{\frac{n-1}{2}} -\frac{t^{n-1}\dot{H}^{n-1}(t)}{H^{n-1}(t)}
 -(n-1)^{\frac{n-1}{2}}\frac{\lambda^{n-1} (1+y^{2})^{n-1} }{(1+\lambda^{2} y^{2})^{n-1}}\\
& \leq \bigg( \frac{t^{2}\dot{H}^{2}(t)}{H^{2}(t)}+
 (n-1)\frac{ (1+y^{2})^{2} }{4y^{2} }  \bigg)^{\frac{n-1}{2}} -(n-1)^{\frac{n-1}{2}}\frac{ (1+y^{2})^{n-1} }{  (2y )^{n-1}} +\frac{t^{n-1}\dot{H}^{n-1}(t)}{H^{n-1}(t)}
\\
 &\leq \frac{n+1}{2}\frac{t^{n-1}\dot{H}^{n-1}(t)}{H^{n-1}(t)} +\frac{(n-1)^{\frac{n-1}{2}}}{2^{n-2}}\frac{ (1+y^{2})^{n-3} }{ y^{n-3}} \frac{t^{2}\dot{H}^{2}(t)}{H^{2}(t)}  .
\end{split}\]
Since the assumption $h^{\lambda}\in \mathcal{P}(\mathbb{A},\mathbb{A_{*}})$ implies that $H\in \mathcal{W}^{1,n-1}(E)$,
we obtain
\[\begin{split}
\int_{r}^{R} \int_{0}^{+\infty}
& \left( \frac{n+1}{2}\frac{t^{n-1}\dot{H}^{n-1}(t)}{H^{n-1}(t)} +\frac{(n-1)^{\frac{n-1}{2}}}{2^{n-2}}\frac{ (1+y^{2})^{n-3} }{ y^{n-3}} \frac{t^{2}\dot{H}^{2}(t)}{H^{2}(t)}
 \right) \frac{y^{n-2}}{(1+y^{2})^{n-1}}dy dt\\
&<\infty.
\end{split}\]
Then it follows from the Lebesgue's dominated convergence theorem that
\[\begin{split}
\lim_{\lambda\rightarrow 0^{+}}&\int_{r}^{R} \int_{0}^{+\infty}
\bigg[\left( \frac{t^{2}\dot{H}^{2}(t)}{H^{2}(t)}+
 (n-1)\frac{\lambda^{2} (1+y^{2})^{2} }{(1+\lambda^{2} y^{2})^{2}} \right)^{\frac{n-1}{2}} -\frac{t^{n-1}\dot{H}^{n-1}(t)}{H^{n-1}(t)} \\
 & - (n-1)^{\frac{n-1}{2}}\frac{\lambda^{n-1} (1+y^{2})^{n-1} }{(1+\lambda^{2} y^{2})^{n-1}}
 \bigg] \frac{y^{n-2}}{(1+y^{2})^{n-1}}dy dt \\
=&\int_{r}^{R} \int_{0}^{+\infty}
\lim_{\lambda\rightarrow 0^{+}}\bigg[\left( \frac{t^{2}\dot{H}^{2}(t)}{H^{2}(t)}+
 (n-1)\frac{\lambda^{2} (1+y^{2})^{2} }{(1+\lambda^{2} y^{2})^{2}} \right)^{\frac{n-1}{2}} -\frac{t^{n-1}\dot{H}^{n-1}(t)}{H^{n-1}(t)} \\
 & - (n-1)^{\frac{n-1}{2}}\frac{\lambda^{n-1} (1+y^{2})^{n-1} }{(1+\lambda^{2} y^{2})^{n-1}}
 \bigg] \frac{y^{n-2}}{(1+y^{2})^{n-1}}dy dt =0.
\end{split}\]
This implies that \eqref{eq-4.9} is true,
and  so, \eqref{eq-4.8} follows.
Combining \eqref{eq-4.3}, \eqref{eq-4.7} and \eqref{eq-4.8}, we see that
Claim \ref{claim-4.1} holds true.

\smallskip

Further, it follows from H\"{o}lder's inequality that
\begin{equation}\label{eq-4.10}
\begin{split}
\int_{r}^{R}\frac{t^{n-1}\dot{H}^{n-1}}{H^{n-1}}dt
 \geq& \left(\int_{r}^{R} \frac{\dot{H}}{H} dt\right)^{n-1}
\left(\int_{r}^{R} \frac{1}{t^{\frac{n-1}{n-2}} } dx\right)^{2-n}\\
 =& \frac{\log^{n-1}\frac{R_{*}}{r_{*}}}{(n-2)^{n-2}}
\frac{Rr}{ (  R^{\frac{1}{n-2} }
-r^{\frac{1}{n-2} })^{n-2}}.
\end{split}
\end{equation}
The equality holds if and only if
$\frac{\dot{H}(t)}{H(t)}
=C_{2,1} t^{-\frac{n-1}{n-2}}$
a.e. in $E$ for some constant $C_{2,1}>0$.
Since $H\in \mathcal{W}^{1,n-1}(E)$ and $H$ is a homeomorphism from $E$  onto $F$.
The same  reasoning as in the discussions of Theorem \ref{thm-1.1} shows that
\begin{equation}\label{eq-4.11}
\qquad\quad
H(t)=r_{*}\left( \frac{R_{*}}{r_{*}}\right)^{\alpha(t)}
\;\;\text{or}\;\;
H(t)=R_{*}\left( \frac{r_{*}}{R_{*}}\right)^{\alpha(t)},
\end{equation}
where $\alpha(t)$ is the mapping from \eqref{eq-1.4}.

\smallskip

Now, we are going to finish the proof of the lemma.
By Claim \ref{claim-4.1}, \eqref{eq-4.8}, \eqref{eq-4.10} and \eqref{eq-4.11}, we see that
\[\begin{split}
 \mathcal{E}[h^{1}]
 >&\omega_{n-1} \int_{r}^{R} \left( (n-1)^{\frac{n-1}{2}} + \frac{t^{n-1}\dot{H}^{n-1}(t)}{H^{n-1}(t)}\right)  dt=\lim_{\lambda\rightarrow 0^{+}}
 \mathcal{E}[h^{\lambda}]\\
 \geq&\omega_{n-1}
\left((n-1)^{\frac{n-1}{2}} (R-r) +\frac{\log^{n-1}\frac{R_{*}}{r_{*}}}{(n-2)^{n-2}}
\frac{Rr}{ (  R^{\frac{1}{n-2} }
-r^{\frac{1}{n-2} })^{n-2}}\right)\\
 =& \lim_{\lambda\rightarrow 0^{+}}
 \mathcal{E}[h_{1}^{\lambda}]
=\lim_{\lambda\rightarrow 0^{+}}
 \mathcal{E}[h_{2}^{\lambda}],
\end{split}\]
where
$h_{1}^{\lambda}$ and
$h_{2}^{\lambda}$ are the mappings from \eqref{eq-1.7}.
The mapping $h_{1}^{\lambda}$ preserves the orientation and $h_{2}^{\lambda}$ changes the orientation.
The proof of the lemma is complete.
\epf

Now, we are going to finish the proof of Theorem \ref{thm-1.2}.

\subsection{Proof of Theorem \ref{thm-1.2} }
Let $h=\rho S\in \mathcal{F}(\mathbb{A},\mathbb{A_{*}})$, where $\rho=|h|$ and $S=\frac{h}{|h|}$.
It follows from \eqref{eq-1.1}$\sim$\eqref{eq-1.3}, Theorem \ref{thm-1.1} and Lemma \ref{lem-4.3} that
\[\begin{split}
 \mathcal{E}[h]
 =&\int_{\mathbb{A}}
\left( \|DS(x)\|^{2}+\frac{|\nabla\rho(x)|^{2}}{ \rho^{2}(x)} \right)^{\frac{n-1}{2}}dx\\
 \geq&\int_{\mathbb{A}} \left( \|D S(x)\|^{n-1}+ \frac{|\nabla\rho(x)|^{n-1}}{\rho^{n-1}(x)} \right)dx
=\mathbb{E}[1,1][h]\\
 \geq&\omega_{n-1}
\left(
(n-1)^{\frac{n-1}{2}}(R-r)+
\frac{1}{(n-2)^{n-2}}  \frac{Rr}{ (  R^{\frac{1}{n-2} }-r^{\frac{1}{n-2} })^{n-2}}\log^{n-1}\frac{R_{*}}{r_{*}}
\right)\\
 =& \lim_{\lambda\rightarrow 0^{+}}
 \mathcal{E}[h_{1}^{\lambda}]
=\lim_{\lambda\rightarrow 0^{+}}
 \mathcal{E}[h_{2}^{\lambda}],
\end{split}\]
where $h_{1}^{\lambda}$
and $h_{2}^{\lambda}$ are the mappings from
 \eqref{eq-1.7}.
  Thus,
 \begin{eqnarray*}
& &\inf_{h\in \mathcal{F}(\mathbb{A},\mathbb{A}_{*})}
 \mathcal{E}[h]
=\inf_{h\in \mathcal{P}(\mathbb{A},\mathbb{A}_{*})}
 \mathcal{E}[h]
=\lim_{\lambda\rightarrow 0^{+}}
 \mathcal{E}[h_{1}^{\lambda}]
=\lim_{\lambda\rightarrow 0^{+}}
 \mathcal{E}[h_{2}^{\lambda}],
\end{eqnarray*}
and the proof of the theorem is complete.
\qed
 \smallskip

\section{Proofs  of Theorems \ref{thm-1.3} and \ref{thm-1.4} }\label{sec-5}
The aim of this section is to prove Theorems \ref{thm-1.3} and \ref{thm-1.4}.

\subsection{Proof  of Theorem \ref{thm-1.3} }
For $\lambda>0$, let $h^{\lambda}(x)=H(t)\Phi^{\lambda}(\xi)\in \mathcal{P}(\mathbb{A},\mathbb{A}_{*})$
and
$h(x)=H(t)\xi\in \mathcal{R}(\mathbb{A},\mathbb{A}_{*})$,
where $x=t\xi$, $t=|x|$ and $\xi\in \mathbb{S}^{n-1}$.
If $\lambda\not=1$,
 we see from \eqref{eq-2.8} that $\|D\Phi^{\lambda}(\xi)\|^{2}$ depends on $\xi$.
Since $\frac{n-1}{2}>1$,
we obtain from \eqref{eq-2.8}, \eqref{eq-4.4} and
the Minkowski inequality (cf. \cite[Page 9]{kuang}) that
\[\begin{split}
\mathcal{E}&[h^{\lambda}]
=\int_{r}^{R} \int_{\mathbb{S}^{n-1}}\left( \frac{t^{2}\dot{H}^{2}(t)}{H^{2}(t)}+ \|D\Phi^{\lambda}(\xi)\|^{2} \right)^{\frac{n-1}{2}}d\sigma(\xi)dt\\
 &< \int_{r}^{R}\left[ \left(\int_{\mathbb{S}^{n-1}} \frac{t^{n-1}\dot{H}^{n-1}(t)}{H^{n-1}(t)}d\sigma(\xi)\right)^{\frac{2}{n-1}}+ \left( \int_{\mathbb{S}^{n-1}} \|D\Phi^{\lambda}(\xi)\|^{n-1} d\sigma(\xi)\right)^{\frac{2}{n-1}}\right]^{\frac{n-1}{2}} dt.
\end{split}\]
Further, it follows from \eqref{eq-2.9} and  \eqref{eq-4.3}  that
\begin{equation}\label{eq-5.1}
\mathcal{E}[h^{\lambda}]
<\omega_{n-1}\int_{r}^{R}\left(n-1+\frac{t^{2}\dot{H}^{2}(t)}{H^{2}(t)}\right)^{\frac{n-1}{2}} dt
=\mathcal{E}[h^{1}]= \mathcal{E}[h ].
\end{equation}
By the arbitrariness of $h$,
we see that the proof of the theorem is complete.
\qed

\medskip

In order to prove Theorem  \ref{thm-1.4}, we shall make some preparation.

Let $x=t\xi \in\mathbb{A}$ and $h(x)=H(t)\xi$,
where $t=|x|$ and $H\in C^{2}(E)\cap \mathcal{W}^{1,n-1}(E)$.
It follows from \eqref{eq-5.1} that
\begin{equation}\label{eq-5.2}
 \mathcal{E}[h ]
 =\omega_{n-1}\int_{r}^{R}\Lambda(t,H,\dot{H}) dt:=\mathcal{H}[H],
\end{equation}
where $\Lambda(t,H,\dot{H})=\left(n-1+\frac{t^{2}\dot{H}^{2}(t)}{H^{2}(t)}\right)^{\frac{n-1}{2}} $.
Then the Euler-Lagrange equation (or equilibrium equation) for the energy integral $ \mathcal{E}[h ]$ is
\begin{equation}\label{eq-5.3}
\frac{d}{dt}\left(\frac{\partial}{\partial\dot{H}}\Lambda\right)=\frac{\partial}{\partial H}\Lambda
\end{equation}
(cf. \cite[Section 1.2]{jost}).
By calculations, \eqref{eq-5.3} reduces to
\begin{equation}\label{eq-5.4}
\begin{split}
 0=&(n-3)\frac{t^{2}\dot{H}}{H^{2}} \left( \frac{t^{2}\dot{H}^{2}}{H^{2}}
+ \frac{t^{3}\dot{H}\ddot{H}}{H^{2}} - \frac{t^{3}\dot{H}^{3}}{H^{3}}  \right)\\
&+\left( \frac{t^{2}\dot{H} }{H^{2}} +n-1  \right)
\cdot\left(  2\frac{t^{2}\dot{H} }{H^{2}}+ \frac{t^{3} \ddot{H}}{H^{2}} - \frac{t^{3}\dot{H}^{2}}{H^{3}} \right).
\end{split}
\end{equation}

\begin{lem}\label{lem-5.1}
Suppose that $n\geq4$,  $H\in C^{2}(E)\cap \mathcal{W}^{1,n-1}(E)$ and $H$ is an increasing mapping from $E$ onto $F$.
Then \eqref{eq-5.4} admits a unique solution
$H_{*} $ satisfying
$H_{*}\in C^{\infty}(E)$ and $\dot{H}_{*}>0$ in $E$.
 Further, we have
\begin{equation}\label{eq-5.5}
\begin{split}
\mathcal{H}[H_{*}]
=&
R\omega_{n-1}\left(w_{*}^{2}(R)+n-1\right)^{\frac{n-1}{2}}
-r\omega_{n-1}\left(w_{*}^{2}(r)+n-1\right)^{\frac{n-1}{2}}\\
&-\tau_{*}(n-1)\omega_{n-1}\big(w_{*}(R)-w_{*}(r)\big),
\end{split}
\end{equation}
where $\tau_{*}$ is the constant from \eqref{eq-5.15} and $w_{*}(t)=t\frac{\dot{H}_{*}(t)}{H_{*}(t)}$ in $E$.
 \end{lem}
\bpf
For $t\in E$ and $H\in C^{2}(E) $,
let $w(t)=t\frac{\dot{H}(t)}{H(t)}$ and $\dot{w}(t)=\frac{dw(t)}{dt}$. Then
$$
\frac{t^{3}\dot{H}\ddot{H}}{H^{2}}=w^{3}-w^{2}+t w \dot{w},
$$
and hence, \eqref{eq-5.4} reduces to
\begin{equation}\label{eq-5.6}
0= w^{3}+(n-1)w+t\dot{w}\big( (n-1) +(n-2)w^{2} \big).
\end{equation}
 Solving the differential equation \eqref{eq-5.6}, we get
\begin{equation}\label{eq-5.7}
 w  (n-1+w^{2})^{\frac{n-3 }{2 }} =\frac{\tau}{ t }  ,
\end{equation}
where $\tau$ is an arbitrary constant.
Since $w=\frac{t\dot{H}}{H}$, $H(r)=r_{*}$ and $H(R)=R_{*}$,
we know that $\tau\neq0$.
Then the fact that $H$ is an increasing mapping in $E$ implies that $\tau>0$, $w>0$ and $\dot{H}>0$ in $E$.

For $w\in \mathbb{R}^{+}$,
let
$$
\phi(w)= w  (n-1+w^{2})^{\frac{n-3 }{2 }}.
$$
For any $\tau_{0}$, $t_{0}>0$,
since $\phi$ is a strictly increasing mapping
in $\mathbb{R}^{+}$,
then there exists a unique point $w_{0}\in\mathbb{R}^{+}$
such that $w_{0}(n-1+w_{0}^{2})^{\frac{n-3}{2}} =\frac{\tau_{0}}{t_{0}} $.
Let
$$
\varphi_{0}(t,\tau,w)
=w(n-1+w^{2})^{\frac{n-3 }{2 }} -\frac{\tau}{ t }.
$$
Obviously, $\varphi_{0}\in C^{1}(\mathbb{R}^{+}\times\mathbb{R}^{+}\times\mathbb{R}^{+})$
and $\frac{\partial}{\partial w}\varphi_{0}(t,\tau,w)>0$.
The implicit function theorem yields that
 there exists a unique  function $\varphi\in C^{1}(\mathbb{R}^{+}\times\mathbb{R}^{+},\mathbb{R}^{+})$ such that
$w=\varphi(t,\tau)$.

 Further, for each fixed $\tau>0$,
it follows from
\eqref{eq-5.7} that
\begin{equation}\label{eq-5.8}
\phi' (w)w'(t)=-\frac{ \tau}{t^{2}}.
\end{equation}
Since $\phi \in C^{\infty}(\mathbb{R}^{+})$ and $w\in C^{1}(\mathbb{R}^{+})$,
we get $(\phi' \circ w) \in C^{1}(\mathbb{R}^{+})$.
This, together with the fact that $\phi'(w)=(n-1+w^{2})^{\frac{n-5}{2}}\big(  n-1+(n-2)w^{2}\big)>0$, implies that $w'\in C^{1}(\mathbb{R}^{+})$.
By differentiating with respect to $t$ on both side of \eqref{eq-5.8}, we get $\phi''(w)\big(w'(t)\big)^{2}+\phi'(w)w''(t)=-\frac{2\tau}{t^{3}}$.
Similarly, we see that $w''\in C^{1}(\mathbb{R}^{+})$.
Continue this process.
Finally, we get that $w\in C^{\infty}(\mathbb{R}^{+})$,
and so, $H\in C^{\infty}(E)$.

On the other hand, \eqref{eq-5.6} is equivalent to
$$
-\frac{\dot{H}}{H }=\frac{ n-1+(n-2)w^{2}}{w^{2}+ n-1 }\dot{w}.
$$
Solving this differential equation, we obtain the general solution
$$
H=\kappa
\exp
\left(-(n-2)w+ (n-3)\sqrt{n-1} \arctan\frac{w}{\sqrt{n-1}}\right),
$$
where $\kappa>0$ is an arbitrary constant
and $w=w(t,\tau)$ is the mapping from \eqref{eq-5.7}.
Hence,
\begin{equation}\label{eq-5.9}
H=\kappa
\exp
\left(-(n-2)w(t,\tau)+ (n-3)\sqrt{n-1} \arctan\frac{w(t,\tau)}{\sqrt{n-1}}\right).
\end{equation}

Since $H(r)=r_{*}$ and $H(R)=R_{*}$,
 we obtain from \eqref{eq-5.9} that
\begin{equation}\label{eq-5.10}
\log\frac{R_{*}}{r_{*}}
=\log\frac{H(R)}{H(r)}
 =\psi(r,R, \tau),
\end{equation}
where
\begin{equation}\label{eq-5.11}
\begin{split}
\psi(r,R, \tau)
=& (n-2)\big(w(r,\tau)-w(R,\tau)\big)\\
&+ (n-3)\sqrt{n-1} \left(\arctan\frac{w(R,\tau)}{\sqrt{n-1}}- \arctan\frac{w(r,\tau)}{\sqrt{n-1}}\right) .
\end{split}
\end{equation}

In the following, we prove that for any $0<r<R<+\infty$ and $0<r_{*}<R_{*}<+\infty$,
there exists a  unique $\tau_{*}=\tau_{*}(r,R,r_{*},R_{*})>0$ such that \eqref{eq-5.10} holds,
where $\tau_{*}=\tau_{*}(r,R,r_{*},R_{*})$ means that the constant $\tau_{*}$
depends only on the quantities $r,R,r_{*}$ and $R_{*}$.

\begin{claim}\label{claim-5.1}
For any $0<r<R<+\infty$ and $\tau>0$,
we have $\frac{\partial}{\partial \tau}\psi(r,R,\tau)>0$.
\end{claim}

It follows from the fact $w\in C^{1}(\mathbb{R}^{+}\times\mathbb{R}^{+})$ that
$\psi\in C^{1}(\mathbb{R}^{+}\times\mathbb{R}^{+}\times\mathbb{R}^{+})$.
Further, we have
\begin{equation}\label{eq-5.12}
\begin{split}\frac{\partial}{\partial \tau}\psi(r,R,\tau)
=&\frac{\partial w(r,\tau)}{\partial \tau} \left((n-2)-\frac{(n-3)(n-1)}{n-1+w^{2}(r,\tau)} \right)\\
&-\frac{\partial w(R,\tau)}{\partial \tau} \left((n-2)-\frac{(n-3)(n-1)}{n-1+w^{2}(R,\tau)} \right).
\end{split}
\end{equation}
By differentiating with respect to $\tau$ on both side of \eqref{eq-5.7},
we get
$$
\frac{\partial w(t,\tau)}{\partial \tau}
= \frac{\big(n-1+ w^{2}\big)(n-1+w^{2})^{-\frac{n-3}{2}}}
{\big(n-1+(n-2)w^{2}\big)t} .
$$
Combining this with \eqref{eq-5.7}, we get
$$
\frac{\partial w(t,\tau)}{\partial \tau}
= \frac{\big(n-1+ w^{2}\big)w}
{\big(n-1+(n-2)w^{2}\big)\tau},
$$
and so,
$$
\frac{\partial w(r,\tau)}{\partial \tau} \left((n-2)-\frac{(n-3)(n-1)}{n-1+w^{2}(r,\tau)} \right)
=\frac{w(r,\tau)}{\tau}.
$$
This and \eqref{eq-5.12} yield that
$$
\frac{\partial}{\partial \tau}\psi(r,R,\tau)
=\frac{w(r,\tau)}{\tau}-\frac{w(R,\tau)}{\tau}.
$$
Since $w>0$, we infer from \eqref{eq-5.6} that $\frac{\partial}{\partial t}w(t,\tau)<0$.
Then
$w(r,\tau)> w(R,\tau)$ and $\frac{\partial}{\partial \tau}\psi(r,R,\tau)>0$.
The claim holds true.

\smallskip

Obviously, \eqref{eq-5.7} implies
$w(t,\tau)= \frac{\tau}{t(n-1)^{\frac{n-3}{2}}}+o(\tau )\rightarrow0^{+}$
as $\tau\rightarrow0^{+}$.
Recall that $\arctan x_{1}-\arctan x_{2}=\arctan\frac{x_{1}-x_{2}}{1+x_{1}x_{2}}$
and $\arctan x=x+o(x)$ as $x\rightarrow0$,
where $x_{1},x_{2}\in \mathbb{R}^{+}$.
Hence, when $\tau\rightarrow 0^{+}$, we infer from \eqref{eq-5.11} that
\begin{equation}\label{eq-5.13}
\psi(r,R, \tau)
= \frac{\tau}{r(n-1)^{\frac{n-3}{2}}}-\frac{\tau}{R(n-1)^{\frac{n-3}{2}}}
+o(\tau ) \rightarrow 0^{+}.
\end{equation}
 When $\tau\rightarrow +\infty$, \eqref{eq-5.7} implies
$w(t,\tau)=(\frac{\tau}{t})^{\frac{1}{n-2}}+o(\tau^{\frac{1}{n-2}})\rightarrow+\infty$,
 and hence, we obtain from \eqref{eq-5.11} that
\begin{equation}\label{eq-5.14}
\psi(r,R, \tau)
=(n-2)\left(\left(\frac{\tau}{r}\right)^{\frac{1}{n-2}}-\left(\frac{\tau}{R}\right)^{\frac{1}{n-2}}\right)
+o(\tau^{\frac{1}{n-2}})
\rightarrow+\infty.
\end{equation}
 Therefore, for any $0<r<R<+\infty$ and $0<r_{*}<R_{*}<+\infty$,
 we infer from Claim \ref{claim-5.1}, \eqref{eq-5.13} and \eqref{eq-5.14} that
there exists a unique  $\tau_{*}=\tau_{*}(r,R,r_{*},R_{*})>0$ such that
\begin{equation}\label{eq-5.15}
\psi(r,R, \tau_{*})=\log\frac{R_{*}}{r_{*}} .
\end{equation}

Further, by letting $t=r$, $H=r_{*}$ and $\tau=\tau_{*}$ in \eqref{eq-5.9},
we can obtain the constant $\kappa_{*}=\kappa_{*}(r,R,r_{*},R_{*})>0$.
Then
$$H_{*}(t)=\kappa_{*}
\exp
\left(-(n-2)w(t,\tau_{*})+ (n-3)\sqrt{n-1} \arctan\frac{w(t,\tau_{*})}{\sqrt{n-1}}\right)$$
is the unique solution to \eqref{eq-5.4} such that
 $H_{*}\in C^{\infty}(E)$ and
 $\dot{H}_{*}>0$ in $E$.

By \eqref{eq-5.2} and the fact that $w=\frac{t \dot{H}}{H}$ and
$w\in C^{\infty}(\mathbb{R}^{+})$, we obtain
\[\begin{split}
\frac{\mathcal{H}[H]}{\omega_{n-1}}
=&\int_{r}^{R}\left(n-1 +w^{2}(t) \right)^{\frac{n-1}{2}}  dt\\
=&t\big(n-1 +w^{2}(t)\big)^{\frac{n-1}{2}}\Big|^{t=R}_{t=r}
 - (n-1)\int_{r}^{R}t\cdot w(t)w'(t)\left(n-1 +w^{2}(t) \right)^{\frac{n-3}{2}} dt.
\end{split}\]
Combining this with \eqref{eq-5.7},
we see that
$$
 \mathcal{H}[H]
=\omega_{n-1}\left(t\big(n-1 +w^{2}(t)\big)^{\frac{n-1}{2}}  -\tau (n-1) w(t)\right)\Big|^{t=R}_{t=r}  .
$$
Then \eqref{eq-5.5} follows from the above equality and \eqref{eq-5.15}.
 The proof of the lemma is complete.
 \epf

\medskip

Now, we are going to prove Theorem \ref{thm-1.4},

\subsection{Proof of Theorem \ref{thm-1.4} }
Let $\mathcal{B}$ denote the family of strictly increasing homeomorphisms from $E $ onto $F $ in $\mathcal{W}^{1, n-1}(E)$.
For any $h \in\mathcal{R}(\mathbb{A},\mathbb{A_{*}})$,
set $h (x)=H(t)\xi$,
where $x=t\xi$ and $t=|x|$.
Then $H\in \mathcal{B}$.
By \eqref{eq-5.1} and \eqref{eq-5.2}, we know that
there exists a mapping $h^{\lambda}(t\xi)=H(t)\Phi^{\lambda}(\xi)\in  \mathcal{P}(\mathbb{A},\mathbb{A_{*}})$ such that
\begin{equation}\label{eq-5.16}
 \mathcal{E}[h^{\lambda}]
 <\omega_{n-1}\int_{r}^{R}\left(n-1+\frac{t^{2}\dot{H}^{2}(t)}{H^{2}(t)}\right)^{\frac{n-1}{2}} dt
 =\mathcal{H}[H]
= \mathcal{E}[h ].
\end{equation}
For every constant $c>0$,  a similar approach as in the proof of \eqref{eq-3.1} gives that
$$
\mathcal{E}\left[\frac{ch}{|h|^2}\right]=\mathcal{E}[h].
$$
Thus, we assume that
$$
\lim_{|x|\rightarrow r}|h(x)|=r_{*}
\;\;\;\text{and}\;\;\;
\lim_{|x|\rightarrow R}|h(x)|=R_{*}.
$$

\medskip

In order to prove the theorem,  we use $\mathcal{X}$ to denote the family of increasing continuous mappings from $E$ onto $F$ in $\mathcal{W}^{1, n-1}(E)$.

\begin{claim}\label{claim-5.2}
There exists $H_{0}\in \mathcal{X}$
such that
$\min_{H\in \mathcal{X}} \mathcal{H}[H]=\mathcal{H}[H_{0}]$.
\end{claim}

 We will divide the proof of the claim into three steps,
 which is based upon the ideas
from \cite[Theorem 3.3]{dac}.
\smallskip

{\bf Step 1.} Compactness.
Let $\{H_{m}\}$ be  a minimizing sequence in $\mathcal{X}$, i.e.,
\begin{equation}\label{eq-5.17}
\lim_{m\rightarrow\infty} \mathcal{H}[H_{m}]
=\mathcal{I}: =\inf_{H\in  \mathcal{X} } \mathcal{H}[H].
\end{equation}
For $t\in E$ and $H\in \mathcal{X}$, we see that there exists a positive constant $C_{3,1}=C_{3,1}(n,E,F)$ such that
$\Lambda$ satisfies the following coercivity condition
\begin{equation}\label{eq-5.18}
\Lambda(t,H,K)
=\left( n-1 +t^{2}\frac{K^{2}}{H^{2}}\right)^{\frac{n-1}{2}}
\geq C_{3,1}|K|^{n-1}.
\end{equation}
Therefore, for $m$ large enough, we obtain from \eqref{eq-5.16}$\sim$\eqref{eq-5.18} that
$$
\mathcal{I}+1\geq \mathcal{H}[H_{m}]=\omega_{n-1}\int_{r}^{R}\Lambda(t,H_{m},\dot{H}_{m})dt\geq  C_{3,1}\omega_{n-1}\int_{r}^{R}|\dot{H}_{m}|^{n-1}dt.
$$
This, together with the uniformly boundedness of $H_{m}$, i.e., $r_{*}\leq H_{m}(t)\leq R_{*}$ in $E$,
implies that there exists a positive constant $C_{3,2}=C_{3,2}(n,E,F)$ such that
$$
\|H_{m}\| _{\mathcal{W}^{1,n-1}(E)}
\leq C_{3,2}.
$$
Since $n-1\geq 3$ and $H_{m}\in \mathcal{X}$, we deduce from
\cite[Exercise 1.4.5]{dac} that
there exists a subsequence (still denoted by $H_{m}$) such that
$H_{m}$ converges weakly to a continuous mapping $H_{0}$ (denoted by $H_{m}\rightharpoonup  H_{0}$)
in $\mathcal{W}^{1,n-1}(E)$ satisfying $H_{0}(r)=r_{*}$ and $H_{0}(R)=R_{*}$ (cf. \cite[Pages 36 and 89]{dac}).

Moreover, it follows from the fact $\|H_{m}\| _{\mathcal{W}^{1,n-1}(E)}\leq C_{3,2}$ and \cite[Pages 161 and  and 175]{jost}
that there exists a subsequence of $H_m$ (still denoted by $H_{m}$)
such that $H_{m}\rightarrow H_0$ a.e. in $E$ as $m\rightarrow\infty$,
which means that there exists a subset $E'\subseteq E$
such that $H_{m}\rightarrow H_0$ in $E \backslash E'$  as $m\rightarrow\infty$, where the measure $m(E')=0$.
For any $x,y\in E \backslash E'$ and $x<y$,
it follows that
$$
 H_{0}(x)=\lim_{m\rightarrow\infty} H_m(x)
\leq \lim_{m\rightarrow\infty} H_m(y)
=H_{0}(y).
$$
This, together with the continuity of $H_{0}$ and the fact $m(E')=0$,
implies that $H_{0}$ is an increasing mapping in $E$.
Therefore, $H_{0}\in \mathcal{X}$.

{\bf Step 2.} Lower semicontinuity.
In the following, we show that $\lim_{m\rightarrow\infty} \mathcal{H}[H_{m}]\geq  \mathcal{H}[H_{0}]$ as
$H_{m}\rightharpoonup H_{0}$ in $\mathcal{W}^{1,n-1}(E)$.
It follows from the fact
\begin{equation}\label{eq-5.19}
\begin{split}
 \partial_{KK} \Lambda (t,H,K)
=&(n-1)(n-3)\frac{t^{4}K^{2}}{H^4}
\left(n-1 +t^{2}\frac{K^{2}}{H^{2}}\right)^{\frac{n-5}{2}}\\
&+(n-1)\frac{t^{2}}{H^2}
\left(n-1 +t^{2}\frac{K^{2}}{H^{2}} \right)^{\frac{n-3}{2}}\\
>&0
\end{split}
\end{equation}
that $\Lambda(t,H,K)$ is strictly convex in $K$.
Thus,
\begin{equation}\label{eq-5.20}
\Lambda(t,H_{m},\dot{H}_{m})>\Lambda(t,H_{m},\dot{H}_0)+\Lambda_{K}(t,H_{m},\dot{H}_0)(\dot{H}_{m}-\dot{H}_0).
\end{equation}

Fix $\delta>0$. Since $H_{m}\rightarrow H_0$ a.e. in $E$  as $m\rightarrow\infty$, then Egorov's theorem asserts that
$H_{m}\rightarrow H_0$ uniformly in $E_{\delta}$  as $m\rightarrow\infty$.
Here $E_{\delta}$ is a measurable subset of $E$ with measure
$m(E-E_{\delta})\leq  \delta$.
Now write
$$
F_{\delta}=\left\{t\in E:   | \dot{H}_0(t)|\leq \frac{1}{\delta}\right\}.
$$
Then $m(E-F_{\delta})\rightarrow0$ as $\delta\rightarrow0^{+}$.
Set $G_{\delta}:=E_{\delta}\cap F_{\delta}$.
Obviously, we have that $m(E-G_{\delta})\rightarrow0$ as $\delta\rightarrow0^{+}$.
Since $H_{0}\in  \mathcal{W}^{1,n-1}(E)$, $H_m(E)=F $ and $H_{m}\rightarrow H_0$ a.e. in $E$,
then the  Lebesgue dominated convergence theorem tells us that
\begin{equation}\label{eq-5.21}
\lim_{m\rightarrow\infty}\int_{G_{\delta}} \Lambda (t,H_{m},\dot{H}_0)dt
=\int_{G_{\delta}}\Lambda (t,H_{0},\dot{H}_0)dt.
\end{equation}

In the following, we will prove
\begin{equation}\label{eq-5.22}
\lim_{m\rightarrow\infty}\int_{G_{\delta}} \Lambda_{K} (t,H_{m},\dot{H}_0) (\dot{H}_m-\dot{H}_0)dt=0.
\end{equation}
For any $\rho\in\overline{ F}=[r_{*},R_{*}]$, it follows from the fact
$ | \dot{H}_0(t)|\leq  \frac{1}{\delta}$ in $G_{\delta}$
that there exists a constant $C_{3,3}=C_{3,3}(\delta,n,E,F)$ such that $\Lambda_{KH}(t,\rho,\dot{H}_{0})\leq C_{3,3} $ in $G_{\delta}$.
Since $H_{m}\rightarrow H_0$ uniformly in $G_{\delta}$ as $m\rightarrow\infty$,
we obtain that $\Lambda_{K}(t,H_{m},\dot{H}_{0}) \rightarrow \Lambda_{K}(t,H_{0},\dot{H}_{0})$ uniformly in $G_{\delta}$ as $m\rightarrow\infty$.
Therefore, for any $\varepsilon>0$, there exists a positive integer $M_{1}=M_{1}(\varepsilon,\delta,n,E,F)$,
such that for any $m>M_{1}$ and $t\in G_{\delta}$,
$$
|\Lambda_{K}(t,H_{m},\dot{H}_{0}) -\Lambda_{K}(t,H_{0},\dot{H}_{0})|<\varepsilon.
$$
Hence, for each $\delta>0$, \eqref{eq-5.22} follows from $\dot{H}_m\rightharpoonup \dot{H}_0$ in $L^{n-1}(E)$ and the fact that
$\Lambda_{K}(t, H_{0}, \dot{H}_0)\in L^{ \frac{n-1}{n -2}}(E)$.

We deduce from \eqref{eq-5.20}$\sim$\eqref{eq-5.22} that
$$
 \lim_{m\rightarrow\infty} \mathcal{H}[H_{m}]
\geq \lim_{m\rightarrow\infty} \omega_{n-1}\int_{G_{\delta}} \Lambda (t,H_{m},\dot{H}_m)dt
\geq\omega_{n-1}\int_{G_{\delta}}\Lambda (t,H_{0},\dot{H}_0)dt.
$$
Let $\delta\rightarrow 0^{+}$.
We recall the monotone convergence theorem to conclude that
$$
 \lim_{m\rightarrow\infty} \mathcal{H}[H_{m}]\geq\omega_{n-1}\int_{r}^{R}\Lambda(t, H_0,\dot{H}_0)dt
= \mathcal{H}[ H_0].
$$

{\bf Step 3.} Since   $\lim_{m\rightarrow\infty} \mathcal{H}[H_{m}]= \inf_{H\in  \mathcal{X} } \mathcal{H}[H]$,
$\lim_{m\rightarrow\infty} \mathcal{H}[H_{m}]
\geq
 \mathcal{H}[ H_0]$
 and $H_{0}\in  \mathcal{X}$,
we deduce that $\mathcal{H}[ H_0]=\min_{H\in  \mathcal{X} } \mathcal{H}[H]$.
The proof of Claim \ref{claim-5.2} is complete.

\medskip

Next, we prove that
\begin{equation}\label{eq-5.23}
\min_{H\in \mathcal{B}}\mathcal{H}[H]=\mathcal{H}[H_{0}].
\end{equation}
By \cite[Proposition 1.2]{rick},  the continuity of $H_0$ and the fact $H_0\in \mathcal{W}^{1,n-1}(E)$,
we know that $H_0$ is absolutely continuous in $E$.
Recall that the fact $H_{0}$ minimizes the functional $\mathcal{H}[H]$ implies that $H_0$ is a solution of
$$ \int_{r}^{R} \left(\frac{\partial}{\partial H}\Lambda(t,H,\dot{H})\cdot \eta
+ \frac{\partial}{\partial\dot{H}}\Lambda(t,H,\dot{H}) \cdot\dot{\eta} \right)dt=0$$
for all $\eta\in AC_0(E,\mathbb{R})$, i.e., $\eta$ is an absolutely continuous mapping from $E$ into $\mathbb{R}$
and we require that there exist $r<r_1\leq R_1<R$
with $\eta(t)=0$ if $t$ is not contained in $[r_1,R_1]$ (cf. \cite[Page 13]{jost}).
Since $\Lambda \in C^{\infty}\big(\mathbb{R}^{+}\times \mathbb{R}^{+}\times (\mathbb{R}^{+}\cup\{0\})\big)$
and $\Lambda_{KK}>0$ in
$\mathbb{R}^{+}\times \mathbb{R}^{+}\times (\mathbb{R}^{+}\cup\{0\})$,
then by \cite[Page 17]{jost} we know that $H_0\in C^{\infty}(E)$.
Hence, \cite[Page 16]{jost} tells us that $H_0$
is a solution of
\eqref{eq-5.3}.
By Lemma \ref{lem-5.1}, we know that  $H_0= H_{*}\in \mathcal{B}$.
This, together with Claim \ref{claim-5.2}, implies that
$$
 \inf_{H\in \mathcal{B}}\mathcal{H}[H]
\geq \min_{H\in \mathcal{X}}\mathcal{H}[H]
=\mathcal{H}[H_{0}]
 \geq \inf_{H\in \mathcal{B}}\mathcal{H}[H].
$$
Therefore, \eqref{eq-5.23} holds true.

\medskip

Now, we are ready to finish the proof of the theorem.
For any $\lambda>0$, let $h_{0}^{\lambda}(t\xi)=H_{0}(t)\Phi^{\lambda}(\xi)$.
For any $\lambda\neq1$, by \eqref{eq-5.1}, \eqref{eq-5.2} and \eqref{eq-5.23}, we obtain that
$$
\min_{h\in \mathcal{R}(\mathbb{A},\mathbb{A_{*}})}\mathcal{E}[h]
=\min_{H\in \mathcal{B}}\mathcal{H}[H]
=\mathcal{H}[H_{0}]
 =\mathcal{E}[h_{0}^{1}]
>
 \mathcal{E}[h_{0}^{\lambda}].
$$
Hence
$$
\inf_{h\in \mathcal{P}(\mathbb{A},\mathbb{A_{*}})}\mathcal{E}[h]
\leq\mathcal{E}[h_{0}^{\lambda}]
<\mathcal{E}[h_{0}^{1}]
=\min_{h\in \mathcal{R}(\mathbb{A},\mathbb{A_{*}})}\mathcal{E}[h]
=\mathcal{H}[H_{*}],
$$
where $\mathcal{H}[H_{*}]$ is the constant from \eqref{eq-5.5}.
The proof of the theorem is complete.
\qed
\smallskip

\vspace*{5mm}
\noindent{\bf Funding.}
 The first author is partially supported by NNSF of China (No. 11671127,  11801166 and 12071121),
 NSF of Hunan Province (No. 2018JJ3327), China Scholarship Council and the construct program of the key discipline in Hunan Province.

\end{document}